\documentclass{SCIS2018}
\usepackage{tikz}
\usepackage{ifpdf,hyperref}
\usepackage{cite}
\usepackage{amsmath,amssymb,amsfonts}
\usepackage{algorithmic}
\usepackage{graphicx}
\usepackage{textcomp}

\def\J{{\bf 1}}

\DeclareMathOperator{\rank}{rank}

\DeclareMathOperator{\Span}{Span}
\DeclareMathOperator{\Col}{Col}
\DeclareMathOperator{\Row}{Row}

\DeclareMathOperator{\lcm}{lcm}

\def\cal{\mathcal}

\def\ra{\rightarrow}

\def\d{\delta}

\def\D{\Delta}

\def\0{{\bf 0}}

\newcommand{\R}{{\mathbb R}}

\def\dsum{\mathop{\sum}\limits}

\newtheorem{thm}{Theorem}[section]
\newtheorem{dfn}[thm]{Definition}
\newtheorem{prp}[thm]{Proposition}
\newtheorem{exa}[thm]{Example}
\newtheorem{lem}[thm]{Lemma}
\newtheorem{cor}[thm]{Corollary}
\newtheorem{rem}[thm]{Remark}

\begin{document}
\ArticleType{RESEARCH PAPER}
\Year{2018}
\Month{January}
\Vol{61}
\No{1}
\DOI{}
\ArtNo{}
\ReceiveDate{}
\ReviseDate{}
\AcceptDate{}
\OnlineDate{}

\title{Design of zero-determinant strategies and its application to networked repeated games}{Design of zero-determinant strategies and its application to networked repeated games}


\author[1]{Daizhan CHENG}{}
\author[2]{Changxi LI}{lichangxi@pku.edu.cn}

\AuthorMark{Cheng D Z}

\AuthorCitation{Cheng D Z, Li C X}



\address[1]{Key Laboratory of Systems and Control, Academy of Mathematics and Systems Sciences,\\ Chinese Academy of Sciences, Beijing {\rm 100190}, China}
\address[2]{Center for Systems and Control, College of Engineering, Peking University, Beijing 100871, China.}

\abstract{Using semi-tensor product (STP) of matrices, the profile evolutionary equation (PEE) for repeated finite games is obtained. By virtue of PEE, the zero-determinant (ZD) strategies are developed for general finite games. A formula is then obtained to design ZD strategies for  general finite games with multi-player and asymmetric strategies. A necessary and sufficient condition is obtained to ensure the availability of the designed ZD strategies. It follows that  player $i$ is able to unilaterally design $k_i-1$ (one less than the number of her strategies) dominating linear relations about the expected payoffs of all players. Finally, the fictitious opponent player is proposed for networked repeated games (NRGs). A technique is proposed to simplify the model by reducing the number of frontier strategies.}


\maketitle


\section{Introduction}
\label{sec:introduction}
In 2012, the zero-determinant (ZD) strategy was firstly proposed by Press and Dyson \cite{pre12}, which shows that in an iterated prisoner's dilemma there exist strategies that dominate any evolutionary opponent. Since then it has attracted considerable attention from game theoretic community as well as computer, information, systems and control  communities. \cite{hao14} calls it ``an underway revolution in game theory", because it reveals that in a repeated game, a player can unilaterally control her opponent's payoff.

A significant development in the following researches is the so-called ``Akin's Lemma" \cite{aki16}, which is a generalization of Press-Dyson's pioneering work without using determinant form.
Akin's original work is about two player two strategy game (prisoner's dilemma). Then various extensions have been done. In \cite{hil15} the ZD strategies of two player two strategy discounted games were discussed. \cite{mca16} considers two player continuous strategy discounted games. A surprising fact is: as a player adopts ZD strategy, her actions restricting to two discrete levels of cooperation are enough to enforce a linear relationship between payoffs of two players even the opponent has infinitely many donation levels to choose. Multiplayer ZD strategies in games with two actions
have been discussed both for undiscounted payoffs \cite{hil14,pan15} and for
discounted payoffs \cite{gov20}. The most general case with multi-player and arbitrary number of strategies is also investigated by {\cite{ued20,tan21}.

Meanwhile, the characteristics of ZD strategies have also be investigated widely. Particularly, the stability of ZD strategy was analyzed in \cite{szo14}; the robustness of ZD strategies has been investigated in \cite{chen14}; the ZD strategies of noisy repeated games were investigated by \cite{hao15}; the influence of misperception on ZD strategies was discussed in \cite{chen22}, and the evolutionary stability of ZD strategies has also been investigated widely \cite{ste13,ada13,hil13}.
ZD strategy technique has also been used for some particular kinds of games, such as  application to public goods games \cite{pan15}, mining pole games \cite{cao19},  snowdrift game \cite{wan17}, etc.

The early works concern more about the ZD strategy design \cite{pre12,hil14,he16}. Most later works focus on general properties of linear relation for average payoffs of players. For instance, \cite{ued20} proved the existence of the solution of a set of linear relations of average payoffs enforced by ZD strategies and the independence of the relations. \cite{tan21} is also concentrated on the existence theorem.

When the applications of ZD strategies to real game problems are concerned, the design formulas or numerical algorithms are necessary. \cite{tah20} provides a method to design ZD strategies for $2\times 2$ asymmetric games. When asymmetric games are considered, the shortage of similar works in previous studies has been claimed in \cite{tah20}.

Another promising development of ZD strategy is its application to networked evolutionary games. For instance, \cite{szo14} considers the extortion strategy under myopic best response arrangement over networks. Several kinds of evolutions of ZD strategy on networked evolutionary games are revealed by \cite{ron15,xu17}.  \cite{tan20} considers the cooperative mining ZD strategy in block-chain networks, etc.

Recently, a new matrix product, called the semi-tensor product (STP) of matrices, was proposed \cite{che11,che12}, and it has been applied to solve some problems in game theory, including the  modeling and analysis of networked repeated games (NRGs) \cite{che15}, providing a formula to verify whether a finite game is potential \cite{che14}, investigating the vector space structure of finite games and its orthogonal decompositions \cite{che17,hao18}, application to traffic congestion games \cite{zha21} and diffusion games \cite{li18}, just to mention a few. Readers who are interested in the STP approach to finite games are referred to a survey paper \cite{chepr}.

Using STP, this paper presents a profile evolutionary equation (PEE) for general finite repeated games, which is essentially the same as the Markov matrix for the memory-one game in \cite{pre12}. Then a detailed design technique  and rigorous proofs are presented for this general case, which are  generalizations  of those proposed firstly by \cite{pre12}. A necessary and sufficient condition is obtained for the availability of ZD strategies. As a by product, we also prove that if a player has $k_i$ strategies she can provide unilaterally $k_i-1$ linear payoff relations using ZD strategies.

Finally, the ZD strategies for networked repeated games (NRGs) are investigated. By proposing and using fictitious opponent player (FOP) a networked repeated game can be transferred to a two player game, where a player, say, player $i$,  plays with the FOP, who represents the whole network except player $i$.  The ZD strategies for player $i$ are designed for $i$ vs FOP.
%
%
%

The rest of this paper is organized as follows: A brief survey on STP is given in Section 2. Then it is used to develop profile evolutionary equation (PEE) of finite repeated games. Finally, some properties of transition matrix of PPE are investigated, which are important for designing ZD strategies. Section 3 deduces a general formula for designing ZD strategies. A necessary and sufficient condition for the designed ZD strategies to be available is presented. Thereafter, some numerical examples are discussed to illustrate the design procedure. The FOP is proposed in Section 4 for NRGs. Using FOP, the technique of ZD strategies becomes applicable to NRGs.
Section 5 is a brief conclusion.

Before ending this section, the notations used in this paper is presented.
${\cal M}_{m\times n}$ is the set of $m\times n$ dimensional real matrices.
$M^*$ is the adjoint matrix of $M$. $\sigma(M)$ is  the set of eigenvalues of $M$.
$\rho(M)$ is   the spectral radius of $M$.
$M>0$ ($M\geq 0$) represents that  all entries of $M$ are positive (non-negative).
$\ltimes$  is  the STP of matrices.
$\Col(A)$ ($\Row(A)$)  is   the set of columns (rows) of ~$A$; $\Col_i(A)$ ($\Row_i(A)$): the $i$-th column (row) of ~$A$.
${\cal D}_k=\{1,2,\cdots,k\}$.
$\d_k^i$  is   the $i$-th column of identity matrix $I_k$. $\d_k^0$ is for a zero vector of dimension $k$.
${\cal B}=\{0,1\}$; and ${\cal B}^k=\{(b_1,\cdots,b_k)^T\;|\;b_i\in {\cal B},\;\forall i\}$.
$\D_k=\Col(I_k)=\left\{\d_k^i\;|\;i=1,\cdots,k\right\}$
$L\in {\cal M}_{m\times n}$ is called a logical matrix, if $\Col(L)\subset \D_m$.
Let $L=\left[\d_m^{i_1},\d_m^{i_2},\cdots,\d_m^{i_n}\right]$, it is briefly denoted by
	$L=\d_m[i_1,i_2,\cdots,i_n]$.
${\cal L}_{m\times n}$ is  the set of $m\times n$ logical matrices.
$\Upsilon^m$ is  the set of $m$ dimensional (column) random vectors. That is, $x=(x_1,x_2,\cdots,x_m)^T\in \Upsilon^m$ means $x_i\geq 0$, $\forall i$, and $\dsum_{i=1}^mx_i=1$.
$\Upsilon_{m\times n}$ is  the set of $m\times n$ (column) random matrices. That is,
$A\in \Upsilon_{m\times n}$, if and only if,  columns of $A$,  i.e., $\Col_j(A)$, $j=1,2,\cdots,n$, are random vectors.
	${\cal G}_{[n;k_1,k_2,\cdots,k_n]}$ is  the set of finite non-cooperative games with $n$ players, and player $i$ has $k_i$ strategies, $i=1,2,\cdots,n$.

\section{Modeling of Finite  Repeated Games}

\subsection{A Brief Survey on STP}

\begin{dfn}\label{d2.1.1}\cite{che11,che12}
	Let $M\in {\cal M}_{m\times n}$, $N\in {\cal M}_{p\times q}$, and
	$t:=\lcm(n,p)$ be the least common multiple of $n$ and $p$.
	Then  the STP   of $M$ and $N$ is defined as
	\begin{align}\label{2.1.1}
		M\ltimes N:=\left( M\otimes I_{t/n}\right)\left(N\otimes I_{t/p}\right)\in {\cal M}_{(mt/n) \times (qt/p)},
	\end{align}
where $\otimes$ is the  Kronecker product.
\end{dfn}

\begin{rem}\label{r2.1.2}
	\begin{itemize}
		\item[(i)] STP is a generalization of conventional matrix product. That is, if $n=p$, then $M\ltimes N=MN$. It is not necessary (and almost impossible) to distinct STP from conventional matrix product,  because in a computing process, the product might shift from one to the other because of the changes of dimensions. Hence in most cases the symbol $\ltimes$ is omitted.
		
		\item[(ii)] As a generalization, STP keeps all major properties of conventional matrix product available, including  associativity, distributivity, etc. All the properties of matrix product used in this paper are the same for both conventional matrix product and STP.
		
		\item[(iii)] Since  conventional matrix product can be considered as a special case of STP,  all the matrix products used in this paper without product symbol are assumed to be STP.
\end{itemize}
\end{rem}

Next, we consider how to express a finite-valued mapping (or logical mapping) into a matrix form using STP.

Let $f:{\cal D}_m\ra {\cal D}_n$ be a mapping from a finite set to another finite set. Then we can identify $j\in {\cal D}_m$ with its vector form $\vec{j}:=\d_{m}^j\in \D_{m}$. In this way, $f$ can be regarded  as a mapping $f: \D_m\ra \D_n$. In the sequel $\vec{j}$ is simply denoted by $j$ again if there is no possible confusion.

\begin{prp}\label{p2.1.3} Let $f:{\cal D}_m\ra {\cal D}_n$. Then there exists a unique matrix $M_f\in {\cal L}_{m\times n}$, called the structure matrix of $f$, such that as the arguments are expressed into their vector forms, we have
	\begin{align}\label{2.1.2}
		f(x)=M_fx.
	\end{align}
\end{prp}

As a corollary,  Proposition \ref{p2.1.3} can be extended into more general form.

\begin{cor}\label{c1.2.4} Let $x_i\in {\cal D}_{k_i}$, $i=1,2,\cdots,n$, $y_j\in {\cal D}_{p_j}$, $j=1,2,\cdots,m$, and $x=\ltimes_{i=1}^nx_i$, $y=\ltimes_{j=1}^my_j$. Assume
	$$
	y_j=f_j(x_1,x_2,\cdots,x_n),\quad j=1,2,\cdots,m,
	$$
	which have their vector forms as
	\begin{align}\label{2.1.4}
		y_j=M_j\ltimes_{i=1}^nx_i, \quad j=1,2,\cdots,m.
	\end{align}
	Then there exists a unique matrix $M_F$, called the structure matrix of the mapping $F=(f_1,\cdots,f_m)$,
	such that
	\begin{align}\label{2.1.5}
		y=M_Fx,
	\end{align}
	where
	$$
	M_F=M_1*M_2*\cdots*M_n \in {\cal L}_{\rho\times \kappa},
	$$
	and $\rho=\prod_{j=1}^mp_j$, $\kappa=\prod_{i=1}^nk_i$,
	and  $*$ is Kratri-Rao product of matrices.
	\footnote{Let $A\in {\cal M}_{s\times n}$, $B\in {\cal M}_{t\times n}$. Then the Khatri-Rao product of $A$ and $B$, denoted by $A*B\in {\cal M}_{st\times n}$, is defined by \cite{che12}
		$$
		\Col_i(A*B)=\Col_i(A)\Col_i(B),\quad i=1,2,\cdots,n.
		$$
	}
	
\end{cor}

Similarly, we have the following result:

\begin{cor}\label{c1.2.5} Let $x_i\in \Upsilon_{k_i}$, $i=1,2,\cdots,n$ and $y_j\in \Upsilon_{p_j}$, $j=1,2,\cdots,m$, and
	\begin{align}\label{2.1.6}
		y_j=M_jx, \quad j=1,2,\cdots,m,
	\end{align}
	where $M_j\in \Upsilon_{p_j\times \kappa}$. If the random variables $y_j, j=1,2,\cdots,m$ are conditional independent on $x_1,x_2,\cdots,x_n.$ Then there exists a unique matrix $M_F$ such that
	\begin{align}\label{2.1.7}
		y=M_Fx,
	\end{align}
	where
	$$
	M_F=M_1*M_2*\cdots*M_n\in \Upsilon_{\rho\times \kappa},
	$$
	which is also called the structure matrix of the mapping $F=(f_1,f_2,\cdots,f_m)$.
\end{cor}

\subsection{PEE of Finite Games}

\begin{dfn}\label{d2.2.1}
	Consider a finite game $G=(N, S,C)$, where
	\begin{itemize}
		\item[(i)] $N=\{1,2,\cdots,n\}$ is the set of players.
		\item [(ii)] $S=\prod_{i=1}^n S_i$ is the profile, where $S_i=\{1,2,\cdots,k_i\},~i=1,2,\cdots,n,$
		is the strategies (or actions) of player $i$.
		\item[(iii)] $C=(c_1,c_2,\cdots,c_n)$, where $c_i:S\ra \R$
		is the payoff (or utility, cost) function of player $i$, $i=1,2,\cdots,n$.
	\end{itemize}
\end{dfn}

The set of such finite games is denoted by ${\cal G}_{[n;k_1,k_2,\cdots,k_n]}$.
A matrix formulation of the repeated game  $G\in {\cal G}_{[n;k_1,k_2,\cdots,k_n]}$ is described as follows \cite{che15}:

\begin{itemize}
	\item[(i)] Identifying $j\in S_i$ with $\d_{k_i}^j\in \D_{k_i}$, then $S_i\sim \D_{k_i}$.
	
	\item[(ii)] Setting $\kappa=\prod_{i=1}^nk_i$, then $S\sim \D_{\kappa}=\prod_{i=1}^n \D_{k_i}$.
	
	\item[(iii)] Let $x_i\in \D_{k_i}$ be the vector form of a strategy for player $i$. Then $x=\ltimes_{i=1}^nx_i\in \D_{\kappa}$ is a profile.
	
	\item[(iv)] For each player's payoff function $c_i$, there exists a unique row vector $V^c_i\in \R^{\kappa}$ such that
	\begin{align}\label{2.2.1}
		c_i(x)=V^c_ix,\quad i=1,2,\cdots,n.
	\end{align}
\end{itemize}

Now consider a repeated game $G^r$ of $G$, which stands for (infinitely) repeated $G$.  Then each player can determine her action at $t+1$ using historic knowledge. It was proved in \cite{pre12} (see also \cite{hao14}) that:`` the shortest memory player sets the rule of the game, which means the long-memory strategies have no advantages over the memory-one strategies". Based on this observation, the strategy updating rule is assumed Markov-like. That is, the strategy of player $i$ at time $t+1$ depends on the profile at $t$ only. Then we have \cite{che15}
\begin{align}\label{2.2.2}
	x_i(t+1)=L_ix(t),\quad i=1,2,\cdots,n.
\end{align}

Two types of strategies are commonly used:
\begin{itemize}
	\item Pure strategy:
	$$
	L_i\in {\cal L}_{k_i\times \kappa},\quad i=1,2,\cdots,n.
	$$
	\item Mixed Strategy:
	$$
	L_i\in \Upsilon_{k_i\times \kappa},\quad i=1,2,\cdots,n.
	$$
\end{itemize}

Multiplying (by STP) all equations in (\ref{2.2.2}) together yields
\begin{align}\label{2.2.3}
	x(t+1)=Lx(t),
\end{align}
where
$$
L=L_1*L_2*\cdots*L_n.
$$

In pure strategy case $L\in {\cal L}_{\kappa\times \kappa}$ and in mixed strategy case $L\in \Upsilon_{\kappa\times \kappa}$.
In mixed strategy case $x(t)$ can be considered as a distribution of profiles at time $t$. If we take into consideration that $\d_{\kappa}^i$ is used to express the $i$-th profile, then $x(t)$ can also be considered as the expected value of profile at time $t$.

In this paper we concern only mixed strategy case. Now what a player can manipulate is his own strategy updating rule. That is, player $i$ can only choose his $L_i$.

We arrange profiles in alphabetic order as
$$
\begin{array}{ccl}
	S&=&\{(s_1,s_2,\cdots,s_n)\;|\;s_i\in S_i,~i=1,2,\cdots,n\}\\
	~&=&\{(1,1,\cdots,1),(1,1,\cdots,2),\cdots,(k_1,k_2,\cdots,k_n)\}\\
	~&:=&\{s^1,s^2,\cdots,s^{\kappa}\}.
\end{array}
$$

Denote the probability of player $i$ choosing strategy $j$ at time $t+1$ under the situation that the profile at time $t$ is $s^r$ as
\begin{align}\label{2.2.4}
	p^r_{i,j}=Prob(x_i(t+1)=j\;|\;x(t)=s^r).
\end{align}
Then we have the strategy evolutionary equation (SEE) of player $i$ as
\begin{align}\label{2.2.5}
	x_i(t+1)=L_ix(t),
\end{align}
where
\begin{align}\label{2.2.6}
	\begin{array}{l}
		L_i=\begin{bmatrix}
			p^1_{i,1}&p^2_{i,1}&\cdots&p^{\kappa}_{i,1}\\
			p^1_{i,2}&p^2_{i,2}&\cdots&p^{\kappa}_{i,2}\\
			\vdots&~&~&~\\
			p^1_{i,k_i}&p^2_{i,k_i}&\cdots&p^{\kappa}_{i,k_i}\\
		\end{bmatrix}\in \Upsilon_{k_i\times \kappa},~i=1,2,\cdots,n.
	\end{array}
\end{align}


According to Corollary \ref{c1.2.5}, we have PEE as
\begin{align}\label{2.2.7}
	x(t+1)=Lx(t),
\end{align}
where the transition matrix
\begin{align}\label{2.2.701}
	L=L_1*L_2*\cdots*L_n.
\end{align}

We give a simple example to calculate $L$.

\begin{exa}\label{e1.2.6} Consider the repeated prisoners' dilemma. Let $p^r_{i,j}$ be the probability of player $i$ taking strategy $j\in \{C,D\}\sim \{1,2\}$ under the condition $s^r$. Then a straightforward computation shows that
	$$
	\begin{array}{l}
    \begin{cases}
		x_1(t+1)=L_1x(t),\\
		x_2(t+1)=L_2x(t),
    \end{cases}
	\end{array}
	$$
	where
	$$
	L_1=\begin{bmatrix}
		p^1_{1,1}&p^2_{1,1}&p^3_{1,1}&p^4_{1,1}\\
		p^1_{1,2}&p^2_{1,2}&p^3_{1,2}&p^4_{1,2}\\
	\end{bmatrix},~
	L_2=\begin{bmatrix}
		p^1_{2,1}&p^2_{2,1}&p^3_{2,1}&p^4_{2,1}\\
		p^1_{2,2}&p^2_{2,2}&p^3_{2,2}&p^4_{2,2}\\
	\end{bmatrix}.
	$$
	Denote by
	$$
	\begin{array}{l}
		p_i=p^i_{1,1},~q_i=p^i_{2,1},\quad i=1,2,3,4.
	\end{array}
	$$
	It follows that
	$$
	\begin{array}{l}
		p^i_{1,2}=1-p^i_{1,1}=1-p_i,~p^i_{2,2}=1-p^i_{2,1}=1-q_i,\quad i=1,2,3,4.
	\end{array}
	$$
	Then we have
	\begin{align}\label{2.2.8}
		\begin{array}{l}
			L=L_1*L_2~~\\
			=\begin{bmatrix}
				p^1_{1,1}p^1_{2,1}&p^2_{1,1}p^2_{2,1}&p^3_{1,1}p^3_{2,1}&p^4_{1,1}p^4_{2,1}\\
				p^1_{1,1}p^1_{2,2}&p^2_{1,1}p^2_{2,2}&p^3_{1,1}p^3_{2,2}&p^4_{1,1}p^4_{2,2}\\
				p^1_{1,2}p^1_{2,1}&p^2_{1,2}p^2_{2,1}&p^3_{1,2}p^3_{2,1}&p^4_{1,2}p^4_{2,1}\\
				p^1_{1,2}p^1_{2,2}&p^2_{1,2}p^2_{2,2}&p^3_{1,2}p^3_{2,2}&p^4_{1,2}p^1_{2,2}\\
			\end{bmatrix}\\
			=\begin{bmatrix}
					p_1q_1&p_1q_2&p_2q_1&p_2q_2\\
					p_1(1-q_1)&p_1(1-q_2)&p_2(1-q_1)&p_2(1-q_2)\\
					(1-p_1)q_1&(1-p_1)q_2&(1-p_2)q_1&(1-p_2)q_2\\
					(1-p_1)(1-q_1)&(1-p_1)(1-q_2)&(1-p_2)(1-q_1)&(1-p_2)(1-q_2)\\
				\end{bmatrix}
		\end{array}
	\end{align}
\end{exa}

\begin{rem}\label{r1.2.7}
	It is easy to verify that the transition matrix in PEE (refer to (\ref{2.2.8})) is essentially the transpose of the Markov matrix for the memory one game in \cite{pre12}. Corresponding to the ``column order" of \cite{pre12} there is a ``row order" change in (\ref{2.2.8}). This is because our profiles are ordered in alphabetic as $CC,~CD, ~DC, ~DD$, while \cite{pre12} uses the order $CC,~DC,~CD,~DD$.
\end{rem}

\subsection{Properties of PEE}

In this subsection we investigate some properties of the transition matrix $L$ of PEE, which are required for designing ZD strategies. As aforementioned in Remark \ref{r1.2.7}, $L$ is the same as the Markov transition matrix for the memory one game in \cite{pre12} (only with a transpose). So $L$ is a column random matrix. This difference does not affect the following discussion. Hence the following argument is a mimic of the corresponding argument in  \cite{pre12}. What we are going to do is to extend it to general case and put it on a solid mathematical foundation.

In the sequel, we need an assumption on $L$. To present it, some preparation is necessary.

A random square matrix $M$ is called a primitive matrix if there exists a finite integer $s>0$ such that $M^s>0$ \cite{hor86}.
Some nice properties of primitive matrix are cited as follows:

\begin{prp}\label{p2.3.1}(Perron-Frobenius Theorem)\cite{hor86} Let $L$ be a primitive stochastic matrix. Then
	\begin{itemize}
		\item[(i)]  $\rho(L)=1$ and there exists a unique $\lambda\in \sigma(L)$ such that $|\lambda|=1$.
		
		\item[(ii)]
		\begin{align}\label{2.3.1}
			\lim_{t\ra \infty}L^t=P>0.
		\end{align}
		Moreover, $P=uv^T$, where $Lu=u$, $u>0$, $L^Tv=v$, $v>0$.
	\end{itemize}
\end{prp}

We are ready to present our fundamental assumption.

{\bf Assumption A-1:} $L$ is primitive.
\begin{rem}\label{r2.3.2}
	\begin{itemize}
		\item[(i)] A-1 is not always true. For instance, consider (\ref{2.2.8}) and let $p^2_{1,1}=0$, $p^3_{1,1}=0$, $p^4_{1,1}=0$, and $p^4_{2,1}=0$. Then $L$ is not primitive.
		\item[(ii)] If $0<p^r_{i,j}<1$, $\forall r,i,j$, then a straightforward verification shows that $L$  is primitive. So A-1 is always true except a zero-measure set.
		\item[(iii)] According to Proposition \ref{p2.3.1}, we have the following immediate conclusions. (a) If $L$ is primitive, then
		\begin{align}\label{2.3.2}
			\rank(L-I_{\kappa})=\kappa-1.
		\end{align}
		(b) There exists $P=uv^T$, where $Lu=u$, $u>0$, $L^Tv=v$, $v>0$, such that (\ref{2.3.1}) holds. That is,
		\begin{align}\label{2.3.3}
			\lim_{t\ra \infty}L^t=uv^T.
		\end{align}
	\end{itemize}
\end{rem}
\begin{prp}\label{p2.3.3} Let $L$ be a $\kappa\times \kappa$ column primitive stochastic matrix.
	Define $M:=L-I_{\kappa}$ and denote by $M^*$  its adjoint matrix. Then
	\begin{itemize}
		\item[(i)]  $\rank(M^*)=1$.
		\item[(ii)]
		\begin{align}\label{2.3.4}
			\Col_j(M^*)\neq 0,\quad j=1,2,\cdots,\kappa.
		\end{align}
	\end{itemize}
\end{prp}

The proof of this proposition and all other proofs can be found in Appendix.
\begin{prp}\label{p2.3.4} Consider the PEE (\ref{2.2.7}). If $L$ is primitive, then
	\begin{align}\label{2.3.5}
		x^*:=\lim_{t\ra \infty}x(t)=u/\|u\|,
	\end{align}
	where $u$ comes from (\ref{2.3.3}).
\end{prp}

Hereafter, we assume $u$ has been normalized. Then $x^*=u$ is the only normalized eigenvector of $L$ corresponding to eigenvalue $1$.
\begin{prp}\label{p2.3.5} Assume $L$ is primitive, then
	\begin{align}\label{2.3.6}
		\Col_j(M^*)\propto u,\quad \forall j.
	\end{align}
\end{prp}

Combining (\ref{2.3.4}) and (\ref{2.3.6}) yields
\begin{align}\label{2.3.7}
	\Col_j(M^*)=\mu_ju,\quad \mu_j\neq 0,\;\; j=1,2,\cdots,\kappa.
\end{align}

\section{Design of ZD Strategies for Repeated Games}

\subsection{A Universal Formula for ZD-strategies}

Consider the transition matrix $L$ of PEE (\ref{2.2.7}).
Recall the finite game $G$. For player $i$ with action $j$, define an indicative vector $\xi_{i,j}\in\R^\kappa$ as follows
\begin{align}\label{1.4}
\xi_{i,j}=\ltimes_{\tau=1}^n\gamma_\tau,
\end{align}
where
$$
\gamma_\tau=
\begin{cases}
\mathbf{1}_{k_\tau},\tau\neq i,\\
\d_{k_i}^j,~\tau=i.\\
\end{cases}
$$
$\xi_{i,j}\in\R^\kappa$ is called  strategy extraction vector, which has the following property.
\begin{lem}\label{lem3vc}
Consider the FRG $G^r$. Strategy extraction vector $\xi_{i,j}\in\R^\kappa$ has the following property
\begin{align}\label{eq24v}
\begin{array}{ccl}
\xi_{i,j}^\top L&=&\sum_{a\in\Phi_{i,j}}\Row_a(L)\\
&=&[p^{1}_{i,j},p^{2}_{i,j},\cdots,p^{\kappa}_{i,j}],~\forall i\in N, \forall j\in A_i.
\end{array}
\end{align}
  where $\Phi_{i,j}=\{a=(a_1,\ldots,a_n)\in A~|~a_i=j\}\subseteq A$.
\end{lem}

\noindent{\it Proof:} According to the definition of $\xi_{i,j}$ and $L$, we have
\begin{align*}
\begin{array}{ccl}
\xi_{i,j}^\top L&=&\xi_{i,j}^\top[\Col_1(L),\Col_2(L),\cdots,\Col_k(L)]\\
&=&[\xi_{i,j}^\top\Col_1(L),\xi_{i,j}^\top\Col_2(L),\cdots,\xi_{i,j}^\top\Col_k(L)],
\end{array}
\end{align*}
where for each column
\begin{align*}
\begin{array}{ccl}
\xi_{i,j}^\top\Col_r(L)&=&(\ltimes_{\tau=1}^n\gamma_\tau^\top)(\ltimes_{s=1}^n\Col_r(L_s))\\
&=&(\otimes_{\tau=1}^n\gamma_\tau^\top)(\otimes_{s=1}^n\Col_r(L_s))\\
&=&\otimes_{s=1}^n(\gamma_s^\top\Col_r(L_s))\\
&=&p^{r}_{i,j}.
\end{array}
\end{align*}
The second equality comes from STP's property.  The third equality comes from the property of Kronecker product.

\hfill $\Box$

\begin{rem}\label{r3.101}
(i) The strategy extraction vector $\xi_{i,j}\in\R^\kappa$ is called ``Repeat" strategy in the existing works \cite{mca16}, \cite{ued20}.
(ii) The purpose of equation (\ref{eq24v}) is to pick out the set of rows from matrix $L$, which involve $p^r_{i,j}$. The row labels of such set are denoted by $\Phi_{i,j}$.  Then $\xi_{i,j}^\top L$ is the summation of the rows in $L$, which are labeled by $\Phi_{i,j}$. For each pair  $(i,j)$, $\xi_{i,j}^\top L$ realizes an elementary (equivalent) transformation for $L$, which results in a row of $L$ which contains $p_{i,j}$ only, i.e., this new row does not involve $p^d_{s,r}$, $(s,r)\neq (i,j)$.
\end{rem}

If $L$ is primitive, then it has a stationary distribution  $\mu\in\Upsilon^\kappa$ satisfying
\begin{align}\label{CV1.4.1}
L\mu=\mu\Leftrightarrow (L-I)\mu=0.
\end{align}
Multiplying $\xi_{i,j}$ to both sides of (\ref{CV1.4.1}) yields that
\begin{align}\label{1.4.2}
[\xi_{i,j}^\top L-\xi_{i,j}^\top]\mu=0.
\end{align}
Let $T_i=[\xi_i^\top L-\xi_i^\top],$ where $\xi_i=[\xi_{i,1},\xi_{i,2},\cdots,\xi_{i,k_i}].$
According to \cite{ued20}, the ZD strategy of player $i$ belong to the interaction of two subspaces.
\begin{dfn}\label{lem.1}\cite{ued20}
The ZD strategy $L_i$ of player $i$ exists if and only if
\begin{align}\label{1.4.4}
\Span(V)\cap\Span(T_i^\top)\neq\{{\bf 0}_\kappa\},
\end{align}
where $V=[{\bf 1}_\kappa,(V_1^c)^\top,(V_2^c)^\top,\cdots,(V_n^c)^\top].$
\end{dfn}

\begin{rem}
Definition \ref{lem.1} can be used to detect whether  a given strategy $L_i$ is a ZD strategy or not. However, it is difficult to design a zero-determinant strategy  for a given game.
\end{rem}

In the following, we only consider how to derive player $i$'s zero-determinant strategy $p_{i,j}$ associated with action $j$ using $\xi_{i,j},$ where $p_{i,j}=[p^1_{i,j},p^2_{i,j},\cdots,p^\kappa_{i,j}].$
A general design formula is presented in the following:

\begin{prp}\label{a3.2} Consider a repeated game $G^r$, where  $G\in {\cal G}_{[n;k_1,k_2,\cdots,k_n]}$. Assume player $i$ is aimed at a set of linear relations on the expected payoffs as
	\begin{align}\label{3.7}
		\begin{array}{l}
			\ell_{i,j}(Ec_1,Ec_2,\cdots,Ec_n,1)=0,
			\quad 1\leq i\leq n; j=1,2,\cdots,k_i-1,
		\end{array}
	\end{align}
	where $\ell_{i,j}$ is a linear function and $Ec_i$ is the expected payoff of player $i$. Then her ZD strategies can be designed as
\begin{align}\label{3.8}
\begin{array}{ccl}	p_{i,j}&=&(p^1_{i,j},p^2_{i,j},\cdots,p^{\kappa}_{i,j})\\		~&=&\mu_{i,j}\ell_{i,j}\left(V^c_1,V^c_2,\cdots,V^c_n,\J_{\kappa}^T\right)+\xi_{i,j},\quad j=1,2,\cdots,k_i-1,
\end{array}
\end{align}
where $\mu_{i,j}\neq 0$ are adjustable parameters.
\end{prp}
\noindent{\it Proof}:
If $p_{i,j}$ satisfies (\ref{3.8}), then we have
\begin{align}\label{3.8.1}
\begin{array}{ccl}	p_{i,j}-\xi_{i,j}^\top&=&\mu_{i,j}\ell_{i,j}\left(V^c_1,V^c_2,\cdots,V^c_n,\J_{\kappa}^T\right)\\ &=&\mu_{i,j}\ell_{i,j}\left(Ec_1,Ec_2,\cdots,Ec_n\right)+c\\
&=&0,
\end{array}
\end{align}
where $c$ is a constant. Equation (\ref{3.8.1}) implies that
$$\ell_{i,j}\left(Ec_1,Ec_2,\cdots,Ec_n\right)=0,~1\leq i\leq n; j=1,2,\cdots,k_i-1.$$
\hfill $\Box$

\begin{rem}
Equation (\ref{3.8}) is a fundamental formula, which provides a convenient way to design ZD strategies for our preassigned purposes.  One may concern the time complexity of the proposed formula.
We point that the complexity is related with the number of players n, the number of strategies for each
player $k_i$, which will be used for the construction of vector $\xi_{i,j}$. To reduce the complexity for designing
ZD strategies, fictitious opponent player method is proposed in Section 4.
\end{rem}
\begin{dfn}\label{d3.8} A set of ZD strategies is permissible, if  the following two conditions are satisfied.
	\begin{itemize}
		\item[(i)]
		\begin{align}\label{3.9}
			0\leq p_{i,j}\leq 1,\quad j=1,2,\cdots,k_i-1.
		\end{align}
		\item[(ii)]
		\begin{align}\label{3.10}
			0\leq \dsum_{j=1}^{k_i-1}p_{i,j}\leq 1.
		\end{align}
	\end{itemize}
\end{dfn}

\begin{rem}\label{r3.9}
	\begin{itemize}
		\item[(i)] It is obvious that permissibility is a fundamental requirement. Non-permissible strategies are meaningless.
		
		\item[(ii)] It is clear that player $i$ can unilaterally design at most $|S_i|-1$ linear relations. Because when
		$p_{i,j}$, $j<|S_i|$ are all determined, $p_{i,|S_i|}$ is uniquely determined by
		$$
		p_{i,|S_i|}=\J^T_{\kappa}-\dsum_{j=1}^{|S_i|-1}p_{i,j}.
		$$
		
\item[(iii)] Indeed, player $i$ can design $|S_i|-1$ linear relations as she wish. This is an advantage of the formula (\ref{3.8}), because it clearly tells how many linear relations a player may design. It was pointed out by \cite{ued20} that ``when the number $M_n$ of possible actions for player $n$ is more than two, player $n$ may be able to employ a ZD strategy with dim $V_n\geq 2$ to simultaneously enforce more than one linear relations. Such a possibility has never been reported in the context of ZD strategies."

		\item[(iv)] Of course, player $i$ needs not to design $|S_i|-1$ relations. If she intends to design $r<|S_i|-1$ relations, equation (\ref{3.10}) has to be modified by reducing the summation to $r$ items.
		
		\item[(v)] The ZD design formula (\ref{3.8}) can be used simultaneously by multi-players, or even all $n$ players.
		\item[(vi)] Equation (\ref{2.3.7}) is extremely important for formula (\ref{3.8}) to be available, because it ensures that each row in the $M=L-I_{\kappa}$ is replaceable by a designed linear relation to get zero determinant.
	\end{itemize}
\end{rem}

Even though a set of ZD strategies is permissible, it may not be available, which means  the goal (\ref{3.7}) may not be reached. We need the following result:

\begin{thm}\label{t3.10} Consider a repeated game, $G^r$, where $G\in {\cal G}_{[n;k_1,k_2,\cdots,k_n]}$. The stationary distribution  exists, if and only if,
	\begin{itemize}
		\item[(i)] there exists a $\mu\in \Upsilon^{\kappa}$ such that
		\begin{align}\label{3.11}
			\lim_{t\ra\infty}L^t=\mu \J^T_{\kappa}.
		\end{align}
		\item[(ii)]
		\begin{align}\label{3.12}
			\rank(L-I_{\kappa})=\kappa-1.
		\end{align}
	\end{itemize}
\end{thm}
\begin{rem}
The existence of stationary distribution $\mu$ is only a sufficient condition for a set
of ZD strategy designed by formula (\ref{3.8}) to be available. As pointed by \cite{aki16}, it can be replaced by $\lim_{t\rightarrow\infty}\sum_{k=1}^tx(k),$  which is the same as $\mu$ provided $\mu$ exists.
\end{rem}

\begin{rem}\label{r3.1001} To see that permissibility is not enough to ensure (\ref{3.11}) and (\ref{3.12}), we recall Example \ref{e1.2.6}. Assume
$$
L_1=\begin{bmatrix}
p^1_{1,1}&0&p^3_{1,1}&0\\
p^1_{1,2}&1&p^3_{1,2}&1\\
\end{bmatrix},~
L_2=\begin{bmatrix}
p^1_{2,1}&0&p^3_{2,1}&0\\
p^1_{2,2}&1&p^3_{2,2}&1\\
\end{bmatrix}.
$$
Then it is easy to verify that
$$
\begin{array}{ccl}
L=L_1*L_2=\begin{bmatrix}
p^1_{1,1}p^1_{2,1}&0&p^3_{1,1}p^3_{2,1}&0\\
p^1_{1,1}p^1_{2,2}&0&p^3_{1,1}p^3_{2,2}&0\\
p^1_{1,2}p^1_{2,1}&0&p^3_{1,2}p^3_{2,1}&0\\
p^1_{1,2}p^1_{2,2}&1&p^3_{1,2}p^3_{2,2}&1\\
\end{bmatrix}
~\sim
\begin{bmatrix}
p^1_{1,1}p^1_{2,1}&p^3_{1,1}p^3_{2,1}&0&0\\
p^1_{1,2}p^1_{2,1}&p^3_{1,2}p^3_{2,1}&0&0\\
p^1_{1,1}p^1_{2,2}&p^3_{1,1}p^3_{2,2}&0&0\\
p^1_{1,2}p^1_{2,2}&p^3_{1,2}p^3_{2,2}&1&1\\
\end{bmatrix}.
\end{array}
$$
where $\sim$ stands for similar, which is caused by swapping second row with third row and second column with third column.
Then it is clear that $L^t$ is always similar to a block lower triangular matrix. Hence (\ref{3.11}) can never be satisfied. While using (\ref{3.8}), by choosing suitable $V^c_i$, $i=1,2$ and parameter $\mu_{1,1}$, a permissible set of ZD strateries can easily be constructed, which provides a counter example to show permissibility is not enough to ensure availability.
\end{rem}

Note that verifying the two conditions in Theorem \ref{t3.10} is not an easy job. Hence we may replace them by the following one.

\begin{cor}\label{c3.11} Consider  an repeated game, $G^r$, where $G\in {\cal G}_{[n;k_1,k_2,\cdots,k_n]}$. Assume the PEE of $G^r$ is (\ref{2.2.3}), where $L$ is  primitive, then the set of ZD strategies designed by formula (\ref{3.8}) is available.
\end{cor}

\begin{rem}\label{r3.12}
	\begin{itemize}
		\item[(i)] Even though  primitivity of $L$ is only a sufficient condition, it is almost necessity because only a zero-measure set of $L$ may not be primitive. That is, if $L$ is not primitive then there must be some
		$(r,i,j)$ with $p^r_{i,j}\in {\cal B}$. So the designer, who intends to use ZD strategies, is better to avoid using such values.
		
		\item[(ii)] Any player can not unilaterally make the conditions in Theorem \ref{t3.10} satisfied. It depends on other players' strategies. What the player $i$ can do is to do her best, that is, to ensure her designed rows, $\xi_{i,j}$, $j=1,2,\cdots,k_i$ are linearly independent. (A Chinese idiom says that ``Mou~Shi~Zai~Ren, Cheng~Shi~Zai~Tian". (Man proposes, God disposes.)) That is the situation for ZD strategy designer.
	\end{itemize}
\end{rem}

\subsection{Numerical Examples}

In the following, we discuss some numerical examples:
\begin{exa}\label{e4.1}
	
	Consider a $G\in {\cal G}_{[3;2,3,2]}$. Since
	$k_1=2$, $k_2=3$, and $k_3=2$, using (\ref{eq24v}), it is easy to calculate that
	\begin{align}\label{4.2}
		\begin{array}{l}
			\Phi_{1,1}=\{1,2,3,4,5,6\},\\
			\Phi_{1,2}=\{7,8,9,10,11,12\},\\
			\Phi_{2,1}=\{1,2,7,8\},\\
			\Phi_{2,2}=\{3,4,9,10\},\\
			\Phi_{2,3}=\{5,6,11,12\},\\
			\Phi_{3,1}=\{1,3,5,7,9,11\},\\
			\Phi_{3,2}=\{2,4,6,8,10,12\}.\\
			\xi_{1,1}=[1,1,1,1,1,1,0,0,0,0,0,0],\\
			\xi_{1,2}=[0,0,0,0,0,0,1,1,1,1,1,1],\\
			\xi_{2,1}=[1,1,0,0,0,0,1,1,0,0,0,0],\\
			\xi_{2,2}=[0,0,1,1,0,0,0,0,1,1,0,0],\\
			\xi_{2,3}=[0,0,0,0,1,1,0,0,0,0,1,1],\\
			\xi_{3,1}=[1,0,1,0,1,0,1,0,1,0,1,0],\\
			\xi_{3,2}=[0,1,0,1,0,1,0,1,0,1,0,1].\\
		\end{array}
	\end{align}
	
	These parameters depend on the type of games, precisely speaking, they depend on the parameters $\{n;k_1,\cdots,k_n\}$ only. They are  independent of particular games.
	
	\begin{itemize}
		\item Pinning Strategy:
		
		Assume the payoff vectors are:
		$$
		\begin{array}{ccl}
			V^c_1&=&[-3,-0.5,6,9,8,7,-4,-4.5,5,6.5,5,7],\\
			V^c_2&=&[4,-1,-5,7.5,2,3.5,8,-4,5,8,9,-2],\\
			V^c_3&=&[9,5,-6,-5.5,5.5,8,8.5,5.5,-0,-3.5,4.5,7].
		\end{array}
		$$
Assume player 2 want to design  pinning strategies which enforce the average payoffs of players 1 and 3 to be
$$
\begin{array}{l}
Ec_1=r_1=4,\\
Ec_2=r_2=-3.
\end{array}
$$
She may choose $\mu_{2,1}=\mu_{2,2}=0.1$ and then		
set
		\begin{align}\label{4.4}
			\begin{array}{ccl}
				p_{2,1}&:=&(0.1)*V^c_1-(0.4)*\J_{12}^T+\xi_{2,1}\\
				~&=&[0.3,0.55,0.2,0.5, 0.4, 0.3, 0.2, 0.15, 0.1,0.25,0.1, 0.3]\\
				p_{2,2}&=&(0.1)*V^c_3+(0.3)*\J_{12}^T+\xi_{2,2}\\
				~&=&[0.6, 0.2, 0.1, 0.15, 0.25, 0.5, 0.55, 0.25,0.7, 0.35, 0.15, 0.4].
			\end{array}
		\end{align}
		
		It is ready to verify that the ZD strategies designed in (\ref{4.4}) are permissible.
		
		\item[(ii)] Extortion Strategy:
		
		Consider a $G\in {\cal G}_{[3;2,3,2]}$ again.  Assume the payoff structure vectors are as follows:
		
		$$
		\begin{array}{ccl}
			V^c_1&=&[16, 11, -4, -8, -2, -10.3,11.4,18.5,1.2,-3,-2.5,1.5],\\
			V^c_2&=&[3, 2, -1, 0, 5, -6, 4, 3, 3, 1, -1, 7],\\
			V^c_3&=&[-2.9,0,6.8,7.1,2,-9.4,-8.2,0.4,4.6,6.1,-2,2.3].
		\end{array}
		$$
		
		Player 2 plans to design an extortion strategy against both players 1 and 3. She may design
		$$
		\begin{array}{l}
			Ec_2-r=k_1(Ec_1-r),\\
			Ec_2-r=k_2(Ec_3-r).
		\end{array}
		$$
		
		To this end, she needs to design
		$$
		\begin{array}{l}
			p_{2,1}-\xi_{2,1}=\mu_1\left[(V^c_2-r\J_{12}^T)-k_1(V^c_1-r\J_{12}^T)\right],\\
			p_{2,3}-\xi_{2,3}=\mu_2\left[(V^c_2-r\J_{12}^T)-k_2(V^c_3-r\J_{12}^T)\right].\\
		\end{array}
		$$
		
		Choosing $\mu_1=0.05$, $\mu_2=0.1$, $r=1$, $k_1=1.1$, $k_2=1.2$,
		it follows that
		\begin{align}\label{4.5}
			\begin{array}{ccl}
				p_{2,1}&=&[0.275,0.5,0.175,0.445,0.365,0.2715,0.178,0.1375,0.0890,0.22,0.0925,0.2725],\\
				p_{2,2}&=&[0.668,0.22,0.104,0.168,0.28,0.548,0.604,0.272,0.768,0.388,0.16,0.444].\\
			\end{array}
		\end{align}
		
		The ZD strategies designed in (\ref{4.5}) are also permissible.
		
	\end{itemize}
	
\end{exa}

\begin{rem}\label{r4.2}
	\begin{itemize}
		\item[(i)] In general, to design a set of permissible ZD strategies is not an easy job. To determine related parameters we need to solve a set of linear inequalities.
		
		\item[(ii)] To verify Lemma \ref{lem3vc}, we calculate the  matrix $M=L-I_{\kappa}$ for Example \ref{e4.1} as follows:
			\begin{align}\label{4.6}
				M=\left[
				\begin{array}{llll}
					p^1_{1,1}p^1_{2,1}p^1_{3,1}-1&p^2_{1,1}p^2_{2,1}p^2_{3,1}&\cdots&p^{12}_{1,1}p^{12}_{2,1}p^{12}_{3,1}\\
					p^1_{1,1}p^1_{2,1}p^1_{3,2}&p^2_{1,1}p^2_{2,1}p^2_{3,2}-1&\cdots&p^{12}_{1,1}p^{12}_{2,1}p^{12}_{3,2}\\
					p^1_{1,1}p^1_{2,2}p^1_{3,1}&p^2_{1,1}p^2_{2,2}p^2_{3,1}&\cdots&p^{12}_{1,1}p^{12}_{2,2}p^{12}_{3,1}\\
					p^1_{1,1}p^1_{2,2}p^1_{3,2}&p^2_{1,1}p^2_{2,2}p^2_{3,2}&\cdots&p^{12}_{1,1}p^{12}_{2,2}p^{12}_{3,2}\\
					p^1_{1,1}p^1_{2,3}p^1_{3,1}&p^2_{1,1}p^2_{2,3}p^2_{3,1}&\cdots&p^{12}_{1,1}p^{12}_{2,3}p^{12}_{3,1}\\
					p^1_{1,1}p^1_{2,3}p^1_{3,2}&p^2_{1,1}p^2_{2,1}p^2_{3,2}&\cdots&p^{12}_{1,1}p^{12}_{2,1}p^{12}_{3,2}\\
					p^1_{1,2}p^1_{2,1}p^1_{3,1}&p^2_{1,2}p^2_{2,1}p^2_{3,1}&\cdots&p^{12}_{1,2}p^{12}_{2,1}p^{12}_{3,1}\\
					p^1_{1,2}p^1_{2,1}p^1_{3,2}&p^2_{1,2}p^2_{2,1}p^2_{3,2}&\cdots&p^{12}_{1,2}p^{12}_{2,1}p^{12}_{3,2}\\
					p^1_{1,2}p^1_{2,2}p^1_{3,1}&p^2_{1,2}p^2_{2,2}p^2_{3,1}&\cdots&p^{12}_{1,2}p^{12}_{2,2}p^{12}_{3,1}\\
					p^1_{1,2}p^1_{2,2}p^1_{3,2}&p^2_{1,2}p^2_{2,2}p^2_{3,2}&\cdots&p^{12}_{1,2}p^{12}_{2,2}p^{12}_{3,2}\\
					p^1_{1,2}p^1_{2,3}p^1_{3,1}&p^2_{1,2}p^2_{2,3}p^2_{3,1}&\cdots&p^{12}_{1,2}p^{12}_{2,3}p^{12}_{3,1}\\
					p^1_{1,2}p^1_{2,3}p^1_{3,2}&p^2_{1,2}p^2_{2,1}p^2_{3,2}&\cdots&p^{12}_{1,2}p^{12}_{2,1}p^{12}_{3,2}-1\\
				\end{array}
				\right]
			\end{align}	
		Then it is easy to verify that $\Phi_{i,j}$ is the rows with each components containing $p^t_{i,j}$ as a factor.
		Moreover, a simple calculation shows that Lemma \ref{lem3vc} is correct.
		
		\item[(iii)] To verify the availability of ZD-strategies in (\ref{4.4}),
		we assume the strategies for player 1 is
		$$
		\begin{array}{ccl}
			p_{1,1}=[0.2,0.3,0.8,0.7,0.5,0.4,0.7,0.9,0.2,0.2,0.1,0.9];
		\end{array}
		$$
		the strategies for player 3 is
		$$
		\begin{array}{ccl}
			p_{3,1}=[0.15,0.2,0.8,0.85,0.2,0.35,0.7,0.9,0.2,0.15,0.55,0.35].
		\end{array}
		$$
		
		Then the strategy profile dynamics is:
		$$
		x(t+1)=Lx(t),\quad t\geq 0,
		$$
		where $L$ is
		\begin{scriptsize}
			$$
			\begin{array}{ccl}
				L&=&\left[
				\begin{array}{llllllllllll}
					0.009&0.033&0.128&0.2975&0.04&0.042&0.098&0.1215&0.004&0.0075&0.0055&0.0945\\
					0.051&0.132&0.032&0.0525&0.16&0.078&0.042&0.0135&0.016&0.0425&0.0045&0.1755\\
					0.018&0.012&0.064&0.0892&0.025&0.07&0.2695&0.2025&0.028&0.0105&0.0083&0.126\\
					0.102&0.048&0.016&0.0158&0.1&0.13&0.1155&0.0225&0.112&0.0595&0.0067&0.234\\
					0.003&0.015&0.448&0.2082&0.035&0.028&0.1225&0.486&0.008&0.012&0.0413&0.0945\\
					0.017&0.06&0.112&0.0367&0.14&0.052&0.0525&0.054&0.032&0.068&0.0338&0.1755\\
					0.036&0.077&0.032&0.1275&0.04&0.063&0.042&0.0135&0.016&0.03&0.0495&0.0105\\
					0.204&0.308&0.008&0.0225&0.16&0.117&0.018&0.0015&0.064&0.17&0.0405&0.0195\\
					0.072&0.028&0.016&0.0383&0.025&0.105&0.1155&0.0225&0.112&0.0420&0.0743&0.014\\
					0.408&0.112&0.004&0.0068&0.1&0.195&0.0495&0.0025&0.448&0.238&0.0607&0.0260\\
					0.012&0.035&0.112&0.0893&0.035&0.042&0.0525&0.054&0.032&0.048&0.3713&0.0105\\
					0.068&0.14&0.028&0.0158&0.14&0.078&0.0225&0.006&0.128&0.272&0.3037&0.0195\\
				\end{array}
				\right].
			\end{array}
			$$
		\end{scriptsize}
		
		It is easy to verify that
		$$
		\rank(L-I_{12})=11.
		$$
		Moreover, we also have that
		$$
		\lim_{t\ra \infty} L^t=u\J^T_{12},
		$$
		where
		$$
		\begin{array}{ccl}
			u&=&[0.0731,0.075,0.0715,0.082,0.126,0.0775,0.0434,\\
			~&~&~0.1002,0.0475, 0.1278,0.0683,0.1077]^T,
		\end{array}
		$$
		which is the normalized eigenvector of  $L$ with respect to its (unique) eigenvalue $1$.
	\end{itemize}
\end{rem}

\section{Application to Networked Repeated Games}

This section considers how to design ZD strategies for a player, $i$, in a networked repeated game (NRG). We propose a method, called a fictitious opponent player (FOP).

\subsection{Fictitious Opponent Player}

\begin{dfn}\label{d5.1.1}\cite{che15} An NRG is a triple $((N,E), G, \Pi)$, where $(N,E)$ is a network graph, where $N$ is the set of players; $G\in {\cal G}_{[2;k,k]}$ is a symmetric game with two players, called the fundamental network game; $\Pi$ is the strategy updating rule, which describes how each player to update his strategies using his neighborhood information.
\end{dfn}

\begin{rem}\label{r5.1.2}
	\begin{itemize}
		\item[(i)] $G\in {\cal G}_{[2;k,k]}$ is symmetric, if $S_1=S_2:=S_0$ and for any $x,y\in S_0$
		$$
		c_1(x,y)=c_2(y,x).
		$$
		\item[(ii)] If $(i,j)\in E$, then players $i$ and $j$ will play game $G$ repeatedly. In this paper only the fixed graph is considered. Since $G$ is symmetric, then  the order of two players does not affect the result.
		
		\item[(iii)] Such an NRG is denoted by $G^{nr}=((N,E), G, \Pi)$.
	\end{itemize}
\end{rem}

Let player $i\in N$, and  $\deg(i)=d$. Then she may consider $N\backslash\{i\}$ as one player, called the FOP of $i$, denoted by $p_{-i}$. Assume $|S_0|=k$, the neighbors' strategies can be considered as the strategies of $p_{-i}$. That is, $p_{-i}$ has totally $k^d$ strategies.

In fact, we need not to distinct different neighbors, hence if $S_0=\{s_1,s_2, \cdots,s_k\}$, then the set of strategies of $p_{-i}$, denoted by $S_{-i}$, is
\begin{align}\label{5.1.1}
	S_{-i}=\{\underbrace{s_1s_1\cdots s_1}_d, \underbrace{s_1s_1\cdots s_2}_d,\cdots,\underbrace{s_ks_k\cdots s_k}_d\}.
\end{align}

Each $s_*\in S_{-i}$ can be expressed as
$$
s_*=(\underbrace{s_1s_1\cdots s_1}_{d_1},\underbrace{s_2s_2\cdots s_2}_{d_2},\cdots,\underbrace{s_ks_k\cdots s_k}_{d_k}),
$$
where $d_i\geq 0$ and $d_1+d_2+\cdots+d_k=d$.
Hence, we can also express $s_*$ by $(d_1,d_2,\cdots,d_k)$, which means $s_j$ has been used by $d_j$ neighbors, $1\leq j\leq k$. Using this notation, we have
\begin{align}\label{5.1.2}
	S_{-i}=\{(d_1,d_2,\cdots,d_k)\;|\; d_j\geq 0,\;\forall j;\;\dsum_{j=1}^kd_j=d\}.
\end{align}

It is easy to verify that define the strategies of $p_{-i}$ in this way, by ignoring the order of neighbors, the total number of strategies is reduced from $k^d$ to
$$
|S_{-i}|=\frac{(k+d-1)!}{(k-1)!d!}.
$$
Hence this treatment reduces the computational complexity.

From the point of view of player $i$, the NRG is equivalent to a game between she and $p_{-i}$, who has the set of strategies $S_{-i}$ defined by (\ref{5.1.2}).
Let $s_*=(d_1,d_2,\cdots,d_k)\in S_{-i}$. Then the payoff functions for $c_i$ and player $p_{-i}$, denoted by $c_{-i}$, are
\begin{align}\label{5.1.3}
	\begin{array}{l}
		c_i(x_i,s_*)=\dsum_{j=1}^kd_jc_i(x_i,s_j),\\
		c_{-i}(x_i,s_*)=\dsum_{j=1}^kd_jc_j(x_i,s_j).\\
	\end{array}
\end{align}

Note that the FOP formulation is particularly suitable for using ZD strategies, because it is not effected by the structure and size of network graph, even though the size might be $\infty$. As long as the stationary distribution of the overall network exists, ZD strategies are still applicable. Moreover, it is easily designable.

\subsection{ZD Strategies for NRGs}

This section considers how to design ZD-strategies for NRGs. We describe the process through two examples.

\begin{exa}
Consider prisoner's dilemma $G$. The two strategies for both players are cooperation ($C$) and defect ($D$). Their payoffs are described in Table \ref{Tab.5.2.1}, where, as a convention, $T>R>P>S$.

\begin{table}[!hbtp]
	\centering \caption{Payoff bi-matrix of prisoner's dilemma}\label{Tab.5.2.1} 
	\begin{tabular}{|c||c|c|}
		\hline $P_1\backslash P_2$&$C$&$D$\\
		\hline
		\hline $C$&$R,~R$&$S,~T$\\
		\hline $D$&$T,~S$&$P,P$\\
		\hline
	\end{tabular}
\end{table}

Consider a networked repeated prisoners' dilemma, denoted by $G^{nr}$. The network graph, depicted by Fig. \ref{Fig.5.2.1}, is non-homogeneous.

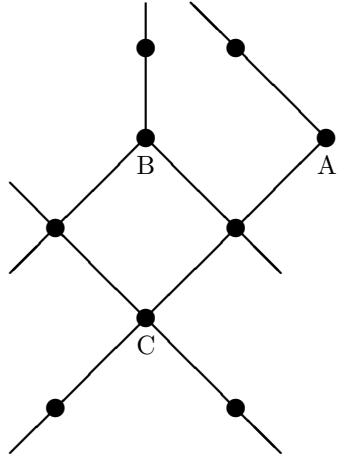
\begin{figure}
	\centering
	\setlength{\unitlength}{6mm}
	\begin{picture}(8,10)\thicklines
		\put(0,0){\line(1,1){7}}
		\put(0,6){\line(1,-1){6}}
		\put(3,7){\line(-1,-1){3}}
		\put(3,7){\line(1,-1){3}}
		\put(3,7){\line(0,1){3}}
		\put(7,7){\line(-1,1){3}}
		\put(1,1){\circle*{0.4}}
		\put(5,1){\circle*{0.4}}
		\put(3,3){\circle*{0.4}}
		\put(1,5){\circle*{0.4}}
		\put(5,5){\circle*{0.4}}
		\put(3,7){\circle*{0.4}}
		\put(7,7){\circle*{0.4}}
		\put(5,9){\circle*{0.4}}
		\put(3,9){\circle*{0.4}}
		\put(2.8,6.2){B}
		\put(6.8,6.2){A}
		\put(2.8,2.2){C}
	\end{picture}
	\caption{A Networked Prisoners' Dilemma}\label{Fig.5.2.1}
\end{figure}

\begin{enumerate}
	
	\item Consider player A. Since $\deg(A)=2$,
	The set of strategies of $p_{-A}$ is
	$$
	S_{-A}=\{(CC),(CD),(DD)\}.
	$$
	Using (\ref{5.1.3}), the payoff vectors for $c_A$ and $c_{-A}$ are, respectively,
	$$
	\begin{array}{l}
		V^c_A=(2R,R+S,2S,2T,T+P,2P),\\
		V^c_{-A}=(2R,R+T,2T,2S,S+P,2P).
	\end{array}
	$$
	It is easy to calculate that $\kappa=6$, and
	$$
	\Phi_{1,1}=\{1,2,3\},~
	\xi_{1,1}=(1,1,1,0,0,0).
	$$
	\begin{itemize}
		\item Pinning Strategy:
		To get $Ec_{-A}=r$, the ZD strategy of player A can be designed as
		$$
		(p^1_{1,1},p^2_{1,1},\cdots,p^{6}_{1,1})=
		\mu(V^c_{-A}-r\J_{6}^T)-\xi_{1,1}.
		$$
		\item Extortion Strategy:
		To get $Ec_A-r=\ell (Ec_{-A}-r)$ with $\ell >1$, the ZD strategy of player A can be designed as
		$$
		\begin{array}{l}
			(p^1_{1,1},p^2_{1,1},\cdots,p^{6}_{1,1})
			=\mu\left( (V^c_{A}-r\J_{6}^T)- \ell (V^c_{-A}-r\J_{6}^T)\right) -\xi_{1,1}.
		\end{array}
		$$
	\end{itemize}
	
	\item Consider player B. Since $\deg(B)=3$,
	The set of strategies of $p_{-B}$ is
	$$
	S_{-B}=\{(CCC),(CCD),(CDD),(DDD)\};
	$$
	Using (\ref{5.1.3}), the payoff vectors for $c_B$ and $c_{-B}$ are, respectively,
	$$
	\begin{array}{ccl}
		V^c_B&=&(3R,2R+S,R+2S,3S,3T,2T+P,T+2P,3P),\\
		V^c_{-B}&=&(3R,2R+T,R+2T,3T,3S,2S+P,S+2P,3P).
	\end{array}
	$$
	
	We  have $\kappa=8$ and
	$$
	\Phi_{1,1}=\{1,2,3,4\},~
	\xi_{1,1}=(1,1,1,1,0,0,0,0).
	$$
	
	The design of ZD strategies is similar to the one for $A$.
	
	\item Consider player C. Since $\deg(C)=4$,
	The set of strategies of $p_{-C}$ is
	$$
	\begin{array}{ccl}
		S_{-C}&=&\{(CCCC),(CCCD),(CCDD),(CDDD),(DDDD)\};
	\end{array}
	$$
	Using (\ref{5.1.3}), the payoff vectors for $c_C$ and $c_{-C}$ are, respectively,
	
	$$
	\begin{array}{ccl}
		V^c_C&=&(4R, 3R+S,2R+2S,R+3S,4S,4T,3T+P,2T+2P,T+3P, 4P),\\
		V^c_{-C}&=&(4R, 3R+T, 2R+2T,R+3T,4T,4S,3S+P,2S+2P,S+3P,4P).
	\end{array}
	$$

	It is easy to calculate that $\kappa=10$ and
	$$
	\Phi_{1,1}=\{1,2,3,4,5\},~
	\xi_{1,1}=(1,1,1,1,1,0,0,0,0,0).
	$$
	
	The design of ZD strategies is similar to the one for $A$, or $B$.
\end{enumerate}
\end{exa}

\textcolor[rgb]{0.00,0.07,1.00}{To illustrate the effectiveness of proposed method for large games, we provide the following example.}
\begin{exa}
\textcolor[rgb]{0.00,0.07,1.00}{Consider a networked repeated prisoners' dilemma, denoted by $G^{nr}$. The network graph, depicted by Fig. \ref{Fig.5.2.2}, is a circular ring with a large  number of nodes.
\begin{figure}[!hbtp]
\centering
\includegraphics[width = 1.6 in]{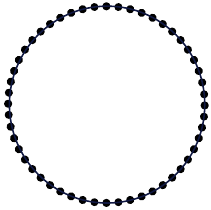}
\caption{The Cycle Ring Graph}
\label{Fig.5.2.2}
\end{figure}}

\textcolor[rgb]{0.00,0.07,1.00}{For any given player $i$, $\deg(i)=2$. The set of strategies of $p_{-i}$ is
$$S_{-i}=\{(CC),(CD),(DD)\}.$$
Using (\ref{5.1.3}), the payoff vectors for $c_i$ and $c_{-i}$ are, respectively,
$$
\begin{array}{l}
V^c_i=(2R,R+S,2S,2T,T+P,2P),~
V^c_{-i}=(2R,R+T,2T,2S,S+P,2P).
\end{array}
$$
It is easy to calculate that $\kappa=6$, and
$\Phi_{1,1}=\{1,2,3\},~\xi_{1,1}=(1,1,1,0,0,0).$}

\textcolor[rgb]{0.00,0.07,1.00}{To realize  $Ec_i-r=0.5(Ec_{-i}-r)$, the ZD strategy of player $i$ can be designed as
		$$
		\begin{array}{l}
			(p^1_{1,1},p^2_{1,1},\cdots,p^{6}_{1,1})
			=\mu\left[ (V^c_{i}-r\J_{6}^T)- 0.5 (V^c_{-i}-r\J_{6}^T)\right] -\xi_{1,1}.
		\end{array}
		$$}
\end{exa}
}

\begin{rem}\label{r5.2.1}
\begin{itemize}
\item[(i)] In the above discussion we only provide the formula for designing ZD strategies. A problem is: is the solution $\{p^i_{j,k},k\in[1,k_i]\}$ obtained from the formula permissible. From designer's point of view, it can be seen immediately from the numerical result. As for theoretical discussion, it is a challenging problem and is out of the scop of this paper. We refer to \cite{ued20} for the existence of the proper solution.

\item[(ii)] When an individual player $i$ in an NRD using ZD strategies, she can ``manipulate" her  immediate neighbors' payoffs
from her. That is the payoffs of the rest of the network got from her. Though it far no means she can manipulate the whole network's payoffs, from her individual perspective, it might be enough.

\item[(iii)] Under our FOP formulation, the ZD strategies in the NRG is exactly the same as  the one for non-networked repeated game.
\item[(iv)]The existence of the zero-determinant strategies is not trivial. We refer to \cite{uedam} for some related
discussion. Further investigation on the existence of ZD strategies for NRGs seems to be necessary and
interesting.
\end{itemize}
\end{rem}

\section{Conclusion}

This paper considers the design of ZD strategies proposed by  Press and  Dyson for general finite games. Using STP, a fundamental formula is presented to numerically realize ZD strategies for finite games with multi-player and asymmetric strategies. In addition to the generality, it simplifies the design procedure.  Then, a necessary and sufficient condition for the availability of the designed ZD strategies is also obtained, which put the ZD technique on a solid foundation. Some numerical examples are presented to demonstrate the efficiency of the method proposed in this paper.

Finally, as an application of the general formula the NRGs are considered. A new concept, called FOP, is proposed as the opponent player for a preassigned player $i$. Using it, the ZD strategies for player $i$ is designed for the game between herself and her FOP. It is interesting that one single player may be able to ``control" the payoff of the rest part of the network from her by using ZD strategies, no matter how large the network is.

\Acknowledgements{The research was supported by National Natural Science Foundation of China under
Grant 62103232, 62073315 and the Natural Science Fund of Shandong Province under grant ZR2021QF005.}

\section*{Appendix}
\begin{enumerate}
	\item

\noindent{\it Proof} of Proposition \ref{p2.3.3}:
		\begin{itemize}
			\item[(i)] Since $L$ has unique eigenvalue $1$, $\rank(M)=\kappa-1$. Hence there exists at least one $(\kappa-1)\times (\kappa-1)$ miner of $M$, which is nonsingular. Hence, $M^*\neq 0$. Observing that
			\begin{align}\label{A.1.1}
				MM^*=\det(M)=0,
			\end{align}
			and $\rank(M)=\kappa-1$, it follows that $\rank(M^*)=1$.
			
			\item[(ii)]  Assume there exists $1\leq j\leq n$ such that $\Col_j(M^*)=0$. Consider $M\backslash \{\Row_j(M)\}$, which is obtained from $M$ by deleting its $j$ th row. Then all its $(\kappa-1)\times (\kappa-1)$ minors have zero determinants. That is, $M\backslash \{\Row_j(M)\}$ is row-dependent.
			
			To get contradiction, we show that any $\kappa-1$ rows of $M$ are linearly independent. Since
			$\dsum_{i=1}^{\kappa}\Row_i(M)={\bf 0}^T_{\kappa}$, $\Row_j(M)=-\dsum_{i\neq j}\Row_i(M)$. If $\rank(M\backslash \{\Row_j(M)\})<\kappa-1$, then $\rank(M)<\kappa-1$, which is a contradiction.
		\end{itemize}
\hfill $\Box$
	
	\item
	
\noindent{\it Proof} of Proposition \ref{p2.3.4}:
		
		First we show that the limit exists. Since $\{x(t)\;|\; t=1,2,\cdots\}\subset \Upsilon^{\kappa}$ and  $\Upsilon^{\kappa}$ is a compact set, there exists a subsequence $\{x_{t_i}\;|\;i=1,2,\cdots\}$ such that
		$$
		\lim_{i\ra \infty}x_{t_i}=x^*\in \Upsilon^{\kappa}.
		$$
		Note that $\lim_{t\ra\infty}L^t=P$, denote by $x^0=Px^*$, we claim that
		\begin{align}\label{A.2.1}
			\lim_{t\ra \infty}x_{t}=x^0.
		\end{align}
		Given any $\epsilon>0$, there exists $N_1$ such that when $t_i>N_1$
		$$
		\|x_{t_i}-x^*\|<\sqrt{\epsilon};
		$$
		and there exists $N_2>0$ such that when $t>N_2$
		$$
		\|M^t-P\|<\sqrt{\epsilon}.
		$$
		Choose  an element $t_{i_0}>N_1$ from the subsequence and set
		$$
		N_3=t_{i_0}>N_1.
		$$
		Assume  $t>N_2+N_3$, then
		$$
		x(t)=M^{t-N_3}x(t_{i_0}).
		$$
		Since $t-N_3>N_2$ and $N_3=t_{i_0}>N_1$, it follows that
		$$
		\|x(t)-x^0\|<(\sqrt{\epsilon})^2=\epsilon.
		$$
		(\ref{A.2.1}) follows. It is also clear that $x^0=x^*$. Moreover, $Lx^*=x^*$ and $Px^*=x^*$. Now since $P=uv^T\in \Upsilon_{\kappa\times \kappa}$, without loss of generality, we can normalize $u$ by replacing $u$ by $u/\|u\|$. Then $v=\J_{\kappa}$. Moreover,
		$$
		x*=Px^*=uv^Tx^*=u.
		$$
		
\hfill $\Box$
	
	\item

\noindent{\it Proof} of Proposition \ref{p2.3.5}:
		
		Note that $MM^*=\det(M)=0$, and $Mu=0$. Since $\rank(M)=\kappa-1$,
		the solutions of equation $Mx=0$ is a one-dimensional subspace. Now each column of $M^*$ is a solution, the conclusion follows.
		
\hfill $\Box$

	\item
	
\noindent{\it Proof} of  Theorem \ref{t3.10}:
		
		(Necessity) It is obvious that $\mu$ is  a stationary distribution if it satisfies the following equation
		\begin{align}\label{A5.1}
			\lim_{t\ra \infty}L^tx_0=u,\quad \forall x_0\in \Upsilon^{\kappa}.
		\end{align}
		We first prove $\lim_{t\ra \infty}L^t$ exists. Since $\Upsilon_{\kappa\times \kappa}$ is a compact set, if
		$\lim_{t\ra \infty}L^t$ does not exist, there must be at least two subsequences $\{L^{n_i}\}$, and $\{L^{m_i}\}$, such that
		$$
		\begin{array}{l}
			\lim_{i\ra \infty}L^{n_i}=P_1,
			\lim_{i\ra \infty}L^{m_i}=P_2,\\
		\end{array}
		$$
		and $P_1\neq P_2$. Say, $\Col_s(P_1)\neq \Col_s(P_2)$, choosing $x_0=\d_{\kappa}^s$, then it violates (\ref{A5.1}).
		
		Hence we have decompositions
		$$
		\lim_{t\ra\infty}L^t=P.
		$$
		Again, because of (\ref{A5.1}) $P$ should have the form that
		$P=[u,u,\cdots,u]$, the conclusion is obvious.
		
		As for the condition (ii), if $\rank(L-I_{\kappa})<\kappa-1$, then $M^*=0$ is a zero matrix.
		Then (\ref{3.7}) fails. Hence (\ref{3.7}) can never be obtained from (\ref{3.8}), and ZD strategies do not work.
		
		(Sufficiency) Replacing any row $s\in \Phi_{i,j}$ of matrix $M=L-I_{\kappa}$ by $\xi_{i,j}$, then condition (ii) ensures (\ref{2.3.7}). Using (\ref{3.8}) and expanding the determinant via replaced row, (\ref{3.7}) follows.
\hfill $\Box$
\end{enumerate}
%


\begin{thebibliography}{99}
%
\bibitem{ada13} Adami C, Hintze A. Evolutionary instability of zero-determinant strategies demonstrates that winning is not everything.
Nature Commun, 2012, 4(1):1-8

\bibitem{aki16} Akin E, The iterated prisoner's delemma: Good strategies and their dynamics. In: Ergodic Theory, Advances in
Dynamical Systems, Berlin, De Gruyter, 2016. 77-107

%
\bibitem{cao19} Cao M, Tang C, Liu Y, Lin F, Chen Z. Application of ZD strategy in mining pool game. In: Proc 38th CCC, 2019. 880-885
%
\bibitem{chen14} Chen J, Zinger A. The robustness of zero-determinant strategies in iterated prisoner's dilemma games. J. Theor.
Biology, 2014, 357: 46-54
\bibitem{che11} Cheng D Z, Qi H S, Li Z Q. Analysis and control of Boolean networks: A semi-tensor product approach. Springer,
London, 2011
\bibitem{che12} Cheng D Z, Qi H S, Zhao Y. An introduction to semi-tensor product of matrices and its applications. World Scientific,
Singapore, 2012
\bibitem{che14} Cheng D Z. On finite potential games. Automatica, 2014, 50(7): 1793-1801
\bibitem{che15} Cheng D Z, He F H, Qi H S, Xu T T. Modeling, analysis and control of networked evolutionary games. IEEE Trans
Aut Contr, 2015, 60(9): 2402-2415
\bibitem{che16} Cheng D Z, Liu T, Zhang K Z. On decomposed subspaces of finite games. IEEE Trans Aut Contr, 2016, 61(11):
3651-3656.
\bibitem{che17} Cheng D Z, Liu T. Linear representation of symmetric games. IET Contr Theory Appl. 2017, 11(18): 3278-3287
\bibitem{chepr} Cheng D Z, Wu Y H, Zhao G D, Fu S H. A comprehensive survey on STP approach to finite games. J Sys Sci Compl,
2021, 34(5): 1666-1680


\bibitem{chen22} Cheng Z, Chen G, Hong Y. Misperception influence on zero-determinant strategies in iterated prisoner's dilemma.
Scientific Reports, 2022, 12(1): 1-9
\bibitem{gov20} Govaert A, Cao M. Zero-determinant strategies in repeated multiplayer social dilemmas with discounted payoffs. IEEE
Trans Aut Contr, 2020, 66(10): 4575-4588
\bibitem{hao18} Hao Y, Cheng D Z. On skew-symmetric games. Journal of the Franklin Institute,2018, 355: 3196-3220
\bibitem{hao14} Hao D, Rong Z H, Zhou T. Zero-determinant strategy: An underway revolution in game theory, Chin Phys B. 2014, 23(7): 078905
\bibitem{hao15} Hao D, Rong Z H, Zhou T. Extortion under uncertainty: zero-determinant strategies in noisy games, Phys Rev E, 2015, 91: 052803
\bibitem{he16} He X, Dai H, Ning P, Dutta R. Zero-determinant strategies for multi-player multi-action iterated games. IEEE Signal
Proc Lett, 2016, 23(3): 311-315
\bibitem{hil13} Hilbe C, Nowak M A, Sigmund K. Evolution of extortion in iterative prisoner's dilemma games. Proc Natl Acad Sci,
2013, 110(17): 6913-6918
\bibitem{hil14} Hilbe C, Wu B, Traulsen A, Nowak M A. Cooperation and control in multplayer social dilemmas. Proc Natl Acad Sci,
2014, 111(46): 16425-16430
\bibitem{hil15} Hilbe C, Traulsen A, Sigmund K. Partners or rivals? Strategies for the iterated prisoner's dilemma. Games and
Economic Behav, 2015, 92: 41-52
\bibitem{hor86} Horn R A, Johnson C R. Matrix Analysis. Cambridge, Campbidge Univ Press, 1986
\bibitem{li18} Li H T, Ding X, Yang Q, Zhou Y. Algebraic formulation and Nash equilibrium of competitive diffusion games. Dyn
Game Appl, 2018, 8: 423-433.
\bibitem{mca16} McAvoy A, Hauert C. Autocratic strategies for iterated games with arbitrary action spaces. Proc Natl Acad Sci, 2016,
113(13): 3573-3578


\bibitem{pan15} Pan L, Hao D, Rong Z, Zhou T. Zero-determinant strategies in iterated public goods game. Science Reports, 2015, 5:
13096
\bibitem{pre12} Press W, Dyson F J. Iterated prosoner's dilemma containes strategies that dominate any evolutionary opponent. Proc
Natl Acad Sci, 2012, 109(26): 10409-10413
\bibitem{ron15} Rong Z, Wu Z, Hao D, et al. Diversity of timescale promotes the maintenance of extortioners in a spatial prisoner's
dilemma game. New Journal of Physics, 2015, 17: 033032.
\bibitem{ste13} Stewart A J, Plotkin J B. Extortion and cooperation in the prisoner's dilemma. Proc Natl Acad Sci, 2013, 110(38):
15348-15353
\bibitem{szo14} Szolnoki A, Perc M. Evolution of extortion in structured populations. Phys Rev E, 2014, 89(2): 022804
\bibitem{tah20} Taha M A, Ghoneim A. Zero-determinant strategies in repeated asymmetric games. Applied Math Comp, 2020, 369: 124862
\bibitem{tan20} Tang C, Li C, Yu X, et al. Cooperative mining in blockchain networks with zero-determinant strategies. IEEE Trans
Cybernetics, 2020, 50(10): 4544-4549
\bibitem{tan21} Tan R, Su Q, Wu B, Wang L. Payoff control in repeated games. In: Proc 33rd Chinese CDC, 2021. 997-1005
\bibitem{ued20} Ueda M, Tanaka T. Linear algebraic structure of zero-determinant strategies in repeated games. PLoS One, 2020,
15(4): e0230973
\bibitem{wan17} Wang J, Guo J, liu H, Shen A. Evolution of zero-determinant strategy on iterated snowdrift game. Acta Phus Sin,
2017, 66(18): 180203
\bibitem{xu17} Xu X, Rong Z, Wu Z, et al. Extortion provides alternative routes to the evolution of cooperation in structured
populations. Physical Review E, 2017, 95(5), 052302.
\bibitem{zha21} Zhang J, Lou J, Qiu J, Lu J. Dynamics and covergence of huper-networked evolutionary games with time delay is
strategies. Information Sciences. 2021, 563: 166-182.
\bibitem{uedam} Ueda M, Necessary and sufficient condition for the existence of zero-determinant strategies in repeated games, Journal
of the Physical Society of Japan 91, 084801, 2022.

%

%

\end{thebibliography}
\end{document}